# The isolated saddle-node bifurcation occurring inside a horseshoe

Yongluo Cao[†] and Shin Kiriki[‡]

[†] Department of Mathematics, Suzhou University,

Suzhou 215006, Jiangsu, P. R. China,

E-mail: `ylcao@suda.edu.cn`

[‡]Department of Mathematical Sciences, Tokyo Denki University,

Hatoyama, Hiki, Saitama, 350-0394, Japan,

E-mail: `ged@r.dendai.ac.jp`

**Abstract**

In this paper, we consider a smooth arc of diffeomorphisms which has a saddle-node bifurcation *inside* a nontrivial invariant set which is a deformation of a horseshoe. We show that this saddle-node bifurcation is *isolated*, that is, its hyperbolicity is maintained before and after the saddle-node bifurcation. Moreover we construct a $C^2$-open set of arcs having the same property.

# 1   Introduction

One can easily imagine such phenomena in nature where some kind of deep impact occurs inside a globally stable system whose stability nevertheless continues to exist after the impact. In this paper, we deal with an example of such phenomena in the horseshoe given by Smale.

This study originates from the consideration of the situation shown in Figure 1. A diffeomorphism with a horseshoe has a configuration of invariant manifolds of a saddle fixed point as shown by the left panel of Figure 1. Some perturbation causes a saddle-node bifurcation in its nontrivial invariant set, that is, a couple of new hyperbolic fixed points is born, and then a new configuration as shown by the right panel of Figure 1 results. In this paper, we clarify a relation between such a saddle-node bifurcation inside the maximal invariant set and the hyperbolicity of its system.

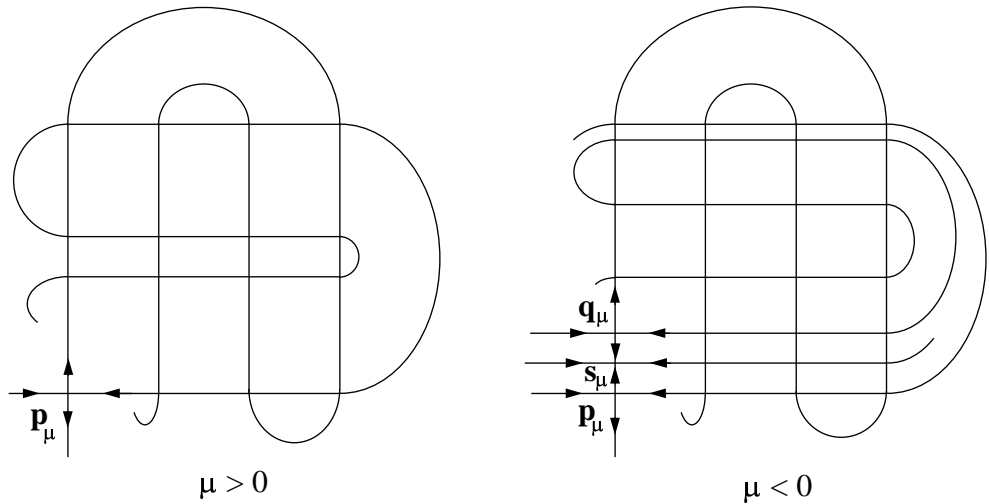

Figure 1

We here consider an arc of diffeomorphisms on a sphere, denoted by $\{\varphi_\mu\}$, such that for $\mu = 1$ it has the nonwandering set which consists of disjoint three basic sets: a hyperbolic nontrivial invariant set called a horseshoe $\Lambda_1$, an attracting fixed point and a repelling fixed point. This is the structurally stable example presented by Smale [8]. Let $\mathbf{p}_1$ be one of two saddle fixed points in this horseshoe such that the stable separatrix $W^s(\mathbf{p}_1)$ and the unstable separatrix $W^u(\mathbf{p}_1)$ are located as shown in the left panel of Figure 1. We assume that this arc has a saddle fixed point $\mathbf{p}_\mu$ for every $\mu$ which is the continuation of $\mathbf{p}_1$. Moreover, we assume that for $\mu = 0$ it has a fixed point $\chi$ such that one of the eigenvalues of $(d\varphi_0)_\chi$ is equal to 1, and the other is smaller than 1. Such a non-hyperbolic fixed point is called an *attracting saddle-node*, and such a parameter value is called a *saddle-node bifurcation point*. After this bifurcation, a new attracting fixed point $\mathbf{s}_\mu$ and a new saddle fixed point $\mathbf{q}_\mu$ appear, as shown in the right panel of Figure 1. That is, before and after this bifurcation, we get a discontinuity of the nonwandering set of $\varphi_\mu$ which is somewhat different from the concepts of explosion or implosion of nonwandering set. If a saddle-node is accumulated by points of a nontrivial invariant set, we say that this saddle-node bifurcation occurs *inside* the nontrivial invariant set. Prof. Robinson dealt with the saddle-node bifurcation occurring *outside* an invariant set which is persistently hyperbolic [7].

The above situation is not special when we consider the elementary bifurcation of 2-dimensional diffeomorphism, because such a phenomenon is obtained from the direct product of some contraction and the well-known arc of one-dimensional maps which has a saddle-node bifurcation as in Figure 2. However, when a parameter value is close to such a bifurcation point, we can not avoid troublesome orbits which need a great many iterations to go through a very small bounded area. That is, they seem to be almost fixed in this small bounded area. Such a characteristic phenomenon caused by the saddle-node bifurcation was first pointed out as an *intermittency* by Pomeau and Manneville, and studied generally by Takens [10]. He classified the generic and attracting saddle-node bifurcations with the intermittent transition into two classes, and moreover suggested the possibility that another significant class exists. In this paper, we consider typical one-parameter families which are within one of the two known classes.

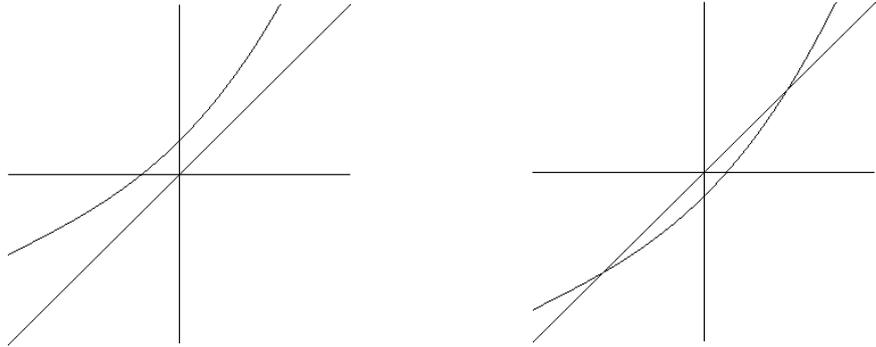

Figure 2

The *first bifurcation* of a parameterized diffeomorphism is an important concept which expresses the first instability of its system, see [6]. It is significant to detect the first bifurcation for given dynamical systems; besides, it is important to know what happens after the first bifurcation. We say that an arc $\{\varphi_\mu\}$ has an *isolated bifurcation* at $\mu = \mu_0$ if there is a neighborhood of $\mu_0$ such that $\varphi_\mu$ is structurally stable for every $\mu$ in this neighborhood except $\mu_0$. When Kiriki studied about arcs of two-dimensional diffeomorphisms having the first homoclinic tangency inside a horseshoe which is a first bifurcation [3], Prof. Takens suggested that the above saddle-node arc may well be an example of the first bifurcation which is different from the first homoclinic tangency [11]. It is clear that such a first homoclinic tangency inside the horseshoe is not an isolated bifurcation.

The aim of this paper is to answer the question: *is the above saddle-node bifurcation inside a nontrivial invariant set a first bifurcation, and moreover, isolated?*

The main result of this paper is

**Main Theorem.** *There exists a $C^2$-open set of smooth arcs of two-dimensional diffeomorphisms such that each of the arcs has a generic and attracting saddle-node bifurcation inside the nontrivial invariant set. Moreover, the bifurcation is isolated.*

We remark that in this paper a boundary cycle of the saddle-node arc in the above class is not critical. However, with some modifications to our arcs, we can make their boundary cycles the *first critical*. The critical of boundary cycle is one of the significant indicators to show that the existence of strange attractors and hyperbolicity are prevalent features, for a generic smooth arcs of diffeomorphisms. See [2].

The above arcs of two-dimensional diffeomorphisms are constructed locally by a direct product of the saddle-node arc of one-dimensional maps and a linear contraction. Therefore, we think that, by the same method, it is possible to construct arcs of diffeomorphisms on a higher dimensional smooth manifold, which have the same property as in the main theorem. For the saddle-node bifurcation occurring in a manifold of dimension greater than or equal to three, there exists an important work by Díaz and Rocha [1].

In section 2, we study the elementary arc of one-dimensional diffeomorphisms. Our two-dimensional saddle-node arcs and a proof of the main theorem is given in section 3.

## 2 Saddle-node arcs of one-dimensional maps

In this section, we consider a one-dimensional saddle-node arc which is essential in the construction of saddle-node arcs of planar diffeomorphisms in the next section. The relationship between the saddle-node arc of one-dimensional maps and the saddle-node vector field has been studied in [9], [5], [2]. Taking these into consideration, we here state clearly its uniformly expanding property near the saddle-node bifurcation.

We first consider an arc of $C^\infty$ diffeomrphisms of $\mathbb{R}$ with coordinate $y$, denoted by $f_\mu(y) = f(\mu, y)$, satisfying

$$f(0, 0) = 0, \quad \frac{\partial}{\partial y}f(0, 0) = 1, \quad \frac{\partial^2}{\partial y^2}f(0, 0) = 2\alpha, \quad \frac{\partial}{\partial \mu}f(0, 0) = \beta,$$

where $\alpha$, $\beta > 0$. For $\mu = 0$, we can get the smooth vector field $Y_0$ defined near $y = 0$ such that $f_0(y) = Y_0(y)_1$ which is a time-one map of $Y_0$, see [9]. In this paper, we assume that, for any positive $\mu$ which is close zero and for any $y$ in a small neighborhood of zero, this arc of diffeomrphisms is embedded as the time-one map of the arc of the saddle-node vector field $Y_\mu(y) = \alpha y^2 + \beta\mu +$ (higher order terms in $\mu$ and $y$). The same assumption is made in [2]. In fact, for any saddle-node arc of diffeomrphisms, Newhouse, Palis and Takens showed the existence of an *adapted* vector field to this arc [5]. We now let $t_1$, $\delta_1$ be sufficiently small constants such that the above saddle-node vector field is defined on $[0, t_1) \times (-\delta_1, \delta_1)$ and moreover, for each $0 \leq \mu < t_1$, there is a point $-\delta_1 < o_\mu < \delta_1$, where $o_0 = 0$, such that $\frac{\partial}{\partial y}Y(\mu, o_\mu) = 0$; $\frac{\partial}{\partial y}Y(\mu, y) > 0$ for any $o_\mu < y < \delta_1$; $\frac{\partial}{\partial y}Y(\mu, y) < 0$ for any $-\delta_1 < y < o_\mu$.

Let $I = [a, b]$ be an interval such that $-\delta_1 < a < o_\mu < b < \delta_1$ and $f_\mu(a) < b$ for every $\mu$ near zero. When $\mu$ is greater than zero, by some iterations of $f_\mu$, the point $a$ can be mapped from the negative side to the positive side. Then, for each $0 < \mu < t_1$, there is a positive integer $n(\mu)$ such that

$$n(\mu) = \min\{n \in \mathbb{Z} \,:\, f_\mu^n(a) \geq b\}.$$

We let the length of the interval $I$ be small so that $f_\mu^{n(\mu)}(b) < \delta_1$. In the next proposition, we show that the derivative of $f_\mu^{n(\mu)}|_I$ uniformly has a positive lower bound.

**Proposition 1** *There exists a positive constant $\sigma_1$ such that for any $a \leq y \leq b$ and any $0 < \mu < t_1$*

$$\frac{d}{dy}f_\mu^{n(\mu)}(y) \geq \sigma_1,$$

*where $n(\mu)$ is defined above.*

**Proof.** For a fixed $0 < \mu < t_1$, let $Y_\mu$ be the above saddle-node field such that $f_\mu$ is embedded as the time-one map of $Y_\mu$. For any $a \leq y \leq b$, we have

$$n(\mu) = \int_y^{f_\mu^{n(\mu)}(y)} \frac{dy}{Y_\mu(y)} = \int_0^{f_\mu^{n(\mu)}(y)} \frac{dy}{Y_\mu(y)} - \int_0^y \frac{dy}{Y_\mu(y)}.$$

Its derivative with respect to $y$ is

$$0 = \frac{1}{Y_\mu(f_\mu^{n(\mu)}(y))} \frac{df_\mu^{n(\mu)}(y)}{dy} - \frac{1}{Y_\mu(y)}.$$

Then, we have

$$\frac{df_\mu^{n(\mu)}(y)}{dy} = \frac{Y_\mu(f_\mu^{n(\mu)}(y))}{Y_\mu(y)}.$$

For each $a \leq y \leq o_\mu$, we have $b \leq f_\mu^{n(\mu)}(a) \leq f_\mu^{n(\mu)}(y) \leq f_\mu^{n(\mu)}(o_\mu)$. Since these points are embedded in the vector field, we have

$$\frac{Y_\mu(f_\mu^{n(\mu)}(y))}{Y_\mu(y)} > \frac{Y_\mu(b)}{Y_\mu(a)}.$$

Since $a, b$ are not singularities of $Y_\mu$ for every $0 < \mu < t_1$, we can define

$$\sigma_0 = \inf \left\{ \frac{Y_\mu(b)}{Y_\mu(a)} \; : \; 0 < \mu < t_1 \right\}.$$

For each $o_\mu < y \leq b$,

$$\frac{Y_\mu(f_\mu^{n(\mu)}(y))}{Y_\mu(y)} \geq \frac{Y_\mu(b)}{Y_\mu(y)} > 1.$$

Then, $\sigma_1 = \min\{\sigma_0, 1\}$ is the desired positive constant. ■

For $-t_1 < \mu < t_1$, we here elucidate the situations of $f_\mu$ outside the interval $(-\delta_1, \delta_1)$. See Figure 3. For each $-t_1 < \mu < t_1$, let $p_\mu$ be a repelling fixed point of $f_\mu$ such that $-\infty < p_\mu < -\delta_1$ and $df_\mu(y)/dy > c_1 > 1$ on a $\delta_2$-neighborhood of $p_\mu$, where $c_1$ and $\delta_2$ are some constants. Except for this repelling fixed point, $f_\mu|_{(-\infty, -\delta_1]}$ has no fixed points for every $-t_1 < \mu < t_1$. Moreover, for each $-t_1 < \mu < t_1$, let $\tilde{p}_\mu$ be a point such that $\delta_1 < \tilde{p}_\mu < +\infty$ and $f_\mu(\tilde{p}_\mu) > \tilde{p}_\mu$, as shown in Figure 3. We assume that, for each $-t_1 < \mu < t_1$ and $p_\mu - \delta_2 \leq y \leq f_\mu(\tilde{p}_\mu)$, $0 < c_2 \leq df_\mu(y)/dy < +\infty$, where $c_2$ is a constant smaller than 1. As before, let $I$ be the same interval $[a, b]$, as in the above construction, which is contained inside the interval $J_\mu = [p_\mu - \delta_2, \ f_\mu(\tilde{p}_\mu)]$. For $0 < \mu < t_1$, $f_\mu|_{J_\mu}$ has no fixed points except the repelling fixed point $p_\mu$ as shown in the left panel of Figure 3.

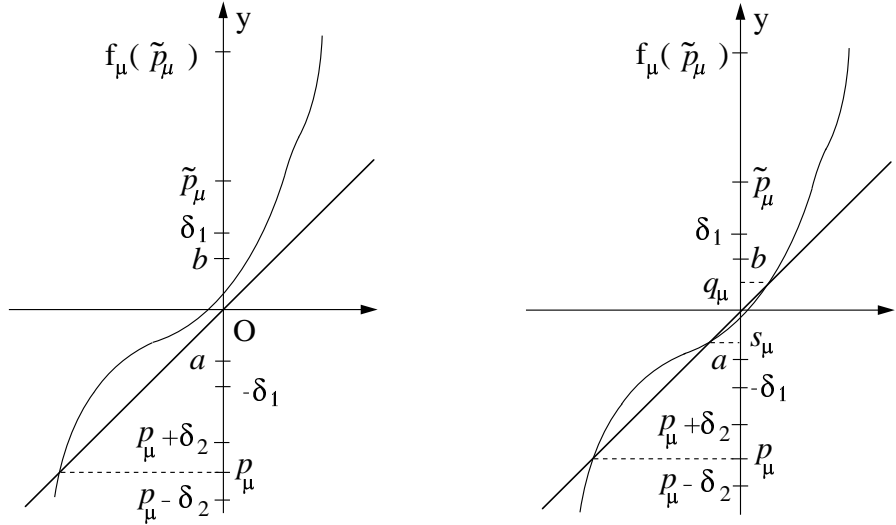

Figure 3

It is clear that, for $\mu < 0$, that is, *after the saddle-node bifurcation*, $f_\mu$ has an attracting fixed point $s_\mu$ and a repelling fixed point $q_\mu$ satisfying $s_\mu < 0 < q_\mu$, which are born as $\mu$ is decreased across the saddle-node bifurcation point $\mu = 0$. Then, there is a constant $t_2$ such that $0 < t_2 \leq t_1$ and , for any $-t_2 < \mu < 0$, both of the hyperbolic fixed points are located inside $I$, as shown in the right panel of Figure 3.

In the following two propositions, we estimate the lower bound for the derivative of some iterations of $f_\mu$ on $J_\mu \setminus I$ for a fixed $-t_2 < \mu < t_1$.

**Proposition 2** . *Let $\mu$ be fixed such that $-t_2 < \mu < t_1$. Then, there exists a constant $\sigma_2 > 0$ such that for any $p_\mu + \delta_2 \leq y < a$ there is a finite iteration number*

$m(y) > 0$ *such that*

$$f_\mu^{m(y)}(y) \geq a; \qquad \frac{d}{dy} f_\mu^{m(y)}(y) \geq \sigma_2.$$

**Proof.** Let $\mu$ be fixed such that $-t_2 < \mu < t_1$. We put

$$r = \inf \{|f_\mu(y) - y| \; : \; p_\mu + \delta_2 \leq y < a\}$$

which is the distance between $f_\mu$ and the identity on $p_\mu + \delta_2 \leq y < a$. We denote the the integer part of $\{a - (p_\mu + \delta_2)\}/r$ by $m$. Moreover, we denote $c_2^{m+1}$ by $\sigma_2$

For any $p_\mu + \delta_2 \leq y < a$, there exists a minimal integer $0 < m(y) < m + 1$ such that

$$f_\mu^{m(y)}(y) \geq a.$$

Since $\left.\frac{df_\mu(y)}{dy}\right|_{f_\mu^i(y)} \geq c_2$ for each $i = 1, \cdots, m(y)$, we have

$$\frac{d}{dy} f_\mu^{m(y)}(y) \geq c_2^{m(y)} > \sigma_2$$

∎

For $-t_2 < \mu < 0$, there exist an attracting fixed point $s_\mu$ and a repelling fixed point $q_\mu$ in the fundamental interval $I = [a, \; b]$. Since $q_\mu$ is repelling and $f_\mu$ is increasing in $I$ for every $-t_1 < \mu < t_1$, we have $df_\mu(y)/dy > 1$ for $q_\mu \leq y \leq b$.

**Proposition 3** . *Let $\mu$ be fixed such that $-t_2 < \mu < t_1$. Then, there exists a constant $\sigma_3 > 0$ such that, for any $b < y \leq f_\mu(\tilde{p}_\mu)$, there is a finite iteration number $\tilde{m}(y) > 0$ such that*

$$f_\mu^{\tilde{m}(y)}(y) > f_\mu(\tilde{p}_\mu); \qquad \frac{d}{dy} f_\mu^{\tilde{m}(y)}(y) \geq \sigma_3.$$

**Proof.** The strategy is essentially the same as that of the proof of the previous proposition. Let $\mu$ be fixed such that $-t_2 < \mu < t_1$. Define

$$\tilde{r} = \inf \{|f_\mu(y) - y| \; : \; b < y \leq f_\mu(\tilde{p}_\mu)\}$$

which is the distance between $f_\mu$ and the identity on $b < y \leq f_\mu(\tilde{p}_\mu)$. We denote the integer part of $(f_\mu(\tilde{p}_\mu) - b)/\tilde{r}$ by $\tilde{m}$, and denote $c_2^{\tilde{m}+1}$ by $\sigma_3$

For any $b < y \leq f_\mu(\tilde{p}_\mu)$, there exists a minimal integer $0 < \tilde{m}(y) < \tilde{m}+1$ such that
$$f_\mu^{\tilde{m}(y)}(y) \geq f_\mu(\tilde{p}_\mu).$$
Since $\left.\frac{df(y)}{dy}\right|_{f_\mu^i(y)} \geq c_2$ for each $i = 1, \cdots, \tilde{m}(y)$, we have
$$\frac{d}{dy} f_\mu^{\tilde{m}(y)}(y) = c_2^{\tilde{m}(y)} \geq \sigma_3$$

∎

We now summarize the properties of the one-dimensional arc $\{f_\mu\}$ *before* the saddle-node bifurcation obtained in the above propositions. For $0 < \mu < t_1$, $f_\mu$ has only the repelling fixed point $p_\mu$, and the derivative of some finite iterations of $f_\mu|_{J_\mu}$ uniformly has a positive lower bound, that is, for each $0 < \mu < t_1$, there are positive integers $n(\mu)$, $m$ and $\tilde{m}$ such that

(i) $\frac{d}{dy} f_\mu(y) > c_1 > 1$ for any $y \in [p_\mu - \delta_2, \ p_\mu + \delta_2)$ since $p_\mu$ is repelling;

(ii) $\frac{d}{dy} f_\mu^{m(y)}(y) > \sigma_2 > 0$ for any $y \in [p_\mu + \delta_2, \ a)$, where $m(y) < m+1$, (from Proposition 2);

(iii) $\frac{d}{dy} f_\mu^{n(\mu)}(y) > \sigma_1 > 0$ for any $y \in [a, \ b]$ (from Proposition 1);

(iv) $\frac{d}{dy} f_\mu^{\tilde{m}(y)}(y) > \sigma_3 > 0$ for any $y \in (b, \ f_\mu(\tilde{p}_\mu)]$, where $\tilde{m}(y) < \tilde{m}+1$, (from Proposition 3).

We next summarize the conditions of the same arc for $-t_2 < \mu < 0$, that is, *after* the saddle-node bifurcation. There exist three hyperbolic fixed points. One is the repelling fixed point $p_\mu$ which is the continuation of the unique repelling fixed point for $\mu > 0$. There is a $\delta_2$-neighborhood of $p_\mu$ such that $df_\mu(y)/dy > 1$ for any $p_\mu - \delta_2 < y < p_\mu + \delta_2$. The others are the repelling fixed point $q_\mu$ and the attracting fixed point $s_\mu$ which are created by the saddle-node bifurcation. There is a neighborhood of $q_\mu$ such that $df_\mu(y)/dy > 1$, and the interval $p_\mu < y < q_\mu$ is a basin of $s_\mu$. By using Proposition 3, for any $q_\mu < y \leq f_\mu(\tilde{p}_\mu)$, we can estimate the norm of the derivative of some iterations of $f_\mu$. Finally, it is clear that the above properties carry on to every arc $C^2$-close enough to $\{f_\mu\}$.

## 3  A saddle-node arc having a horseshoe

Let $\{f_\mu\}$ be the same generic saddle-node arc of one-dimensional $C^\infty$ diffeomorphisms. In this section, using this arc, we construct a generic and attracting saddle-node arc $\{\varphi_\mu\}$ of $C^\infty$ diffeomorphisms on $\mathbb{R}^2$. This arc passes through a somewhat unusual horseshoe map, see Figure 4, so that its expanding factors are given some devices to ensure its hyperbolicity near the saddle-node bifurcation point. In fact, it has an extraordinarily expanding factor in a paticular region. We can construct this arc as follows. The constants used to charactrize $f_\mu$ in the previous section will be used in this section without restating their definitions. $f_\mu$ is again defined on the interval $J_\mu = [p_\mu - \delta_2,\ f_\mu(\tilde{p}_\mu)]$ for $-t_2 < \mu < t_1$. We now consider the rectangle $R_\mu$ in $\mathbb{R}^2$ and the horizontal strip $H_\mu \subset R_\mu$ for $-t_2 < \mu < t_1$ defined as follows:

$$R_\mu = \left\{(x,\ y) \in \mathbb{R}^2\ :\ p_\mu - \delta_2 \leq x \leq f_\mu(\tilde{p}_\mu),\ p_\mu - \delta_2 \leq y \leq f_\mu(\tilde{p}_\mu)\right\};$$

$$H_\mu = \left\{(x,\ y) \in R_\mu\ :\ p_\mu - \delta_2 \leq x \leq f_\mu(\tilde{p}_\mu),\ f_\mu^{-1}(p_\mu - \delta_2) \leq y \leq \tilde{p}_\mu\right\}.$$

Let $\lambda > 0$ be a sufficiently small real number. Its exact range will be given later. For each $-t_2 < \mu < t_1$ and $\mathbf{x} = (x,\ y) \in H_\mu$, we let

$$\varphi(\mu,\ \mathbf{x}) = \varphi_\mu(\mathbf{x}) = \left(\lambda(x - p_\mu) + p_\mu,\ f_\mu(y)\right).$$

Then, we have

$$(d\varphi_\mu)|_{H_\mu} = \begin{pmatrix} \lambda & 0 \\ 0 & \frac{d}{dy}f_\mu(y) \end{pmatrix}.$$

Since $df_\mu(p_\mu)/dy > c_1 > 1$, the point $\boldsymbol{p}_\mu = (p_\mu,\ p_\mu) \in H_\mu$ is a saddle fixed point of $\varphi_\mu$ for every $-t_2 < \mu < t_1$. $\varphi_\mu$ maps the horizontal rectangle $H_\mu$ to a vertical rectangle $\varphi_\mu(H_\mu)$ in $R_\mu$, shown as Figure 4. The bifurcation type at the origin is inherited from the previous one-dimensional saddle-node arc, that is, $\{\varphi_\mu\}$ has a generic and attracting saddle-node bifurcation at $(x,\ y) = (0,\ 0)$ when $\mu = 0$.

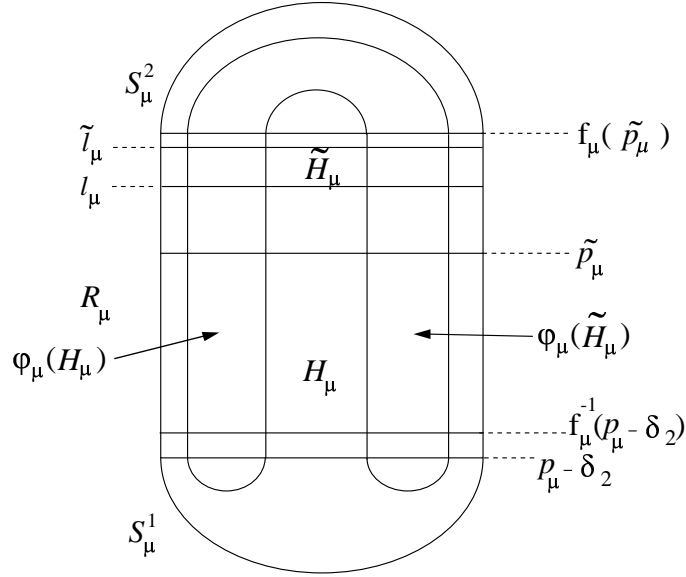

Figure 4

Let $\tilde{\sigma}$ be a sufficiently large real number satisfying the following two conditions: there exists a real constant $\zeta$ such that

$$\tilde{\sigma} \cdot \sigma_1 \cdot \sigma_2 \cdot \sigma_3 \geq \zeta > 1,$$

where $\sigma_1$, $\sigma_2$ and $\sigma_3$ are defined in the previous section, and for every $-t_2 < \mu < t_1$

$$\frac{f_\mu(\tilde{p}_\mu) - (p_\mu - \delta_2)}{\tilde{\sigma}} < f_\mu(\tilde{p}_\mu) - \tilde{p}_\mu.$$

The second condition pronpts us to define the horizontal strip

$$\widetilde{H}_\mu = \left\{(x,\ y) \in R_\mu\ :\ p_\mu - \delta_2 \leq x \leq f_\mu(\tilde{p}_\mu),\ l_\mu \leq y \leq \tilde{l}_\mu\right\},$$

where $l_\mu$ and $\tilde{l}_\mu$ are positive real numbers satisfying

$$\tilde{p}_\mu < l_\mu < \tilde{l}_\mu < f_\mu(\tilde{p}_\mu);$$

$$\tilde{l}_\mu - l_\mu = \frac{f_\mu(\tilde{p}_\mu) - (p_\mu - \delta_2)}{\tilde{\sigma}},$$

as shown in Figure 4. On $\widetilde{H}_\mu$, $\varphi_\mu$ is defined for any $-t_2 < \mu < t_1$ as a linear map satisfying

$$(d\varphi_\mu)|_{\widetilde{H}_\mu} = \begin{pmatrix} -\lambda & 0 \\ 0 & -\tilde{\sigma} \end{pmatrix}$$

and

$$\pi_y \circ \varphi_\mu(x, \, l_\mu) = f_\mu(\tilde{p}_\mu), \quad \pi_y \circ \varphi_\mu(x, \, \tilde{l}_\mu) = p_\mu - \delta_2,$$

where $\pi_y$ is the projection to the vertical axis. Here, we let the above contracting factor $\lambda$ be a small constant satisfying

$$0 < \lambda < \min\left\{\frac{c_2}{2}, \, \frac{1}{\zeta}\right\},$$

where $c_2 = \inf\{df_\mu(y)/y \, : \, -t_2 < \mu < t_1, \, y \in J\}$ as defined in the previous section.

Let $S_\mu^1$ and $S_\mu^2$ be the semidiscs which adhere to the lower side and the upper side of $R_\mu$, respectively, and we assume that $\varphi_\mu(S_\mu^1 \cup S_\mu^2) \subset S_\mu^1$ as ahown in Figure 4. Moreover, the strip between $H_\mu$ and $\widetilde{H}_\mu$ is mapped by $\varphi_\mu$ to the arch connecting $\varphi_\mu(H_\mu)$ and $\varphi_\mu(\widetilde{H}_\mu)$ which is contained in $S_\mu^2$. Now the construction of $\varphi_\mu$ on the rectangle $R_\mu$ complete and $\varphi_\mu(R_\mu)$ is of a horseshoe-shape.

Next, we study the stability of $\varphi_\mu$ when it is near the saddle-node bifurcation. For $0 < \mu < t_1$, that is, before the saddle-node bifurcation, we define a maximal invariant set for $\varphi_\mu$ as

$$\Lambda_\mu = \bigcap_{n \in \mathbb{Z}} \varphi_\mu^n(R_\mu).$$

For each $\mathbf{x} \in \Lambda_\mu$, let $E^s(\mathbf{x})$ be the one-dimensional subspace of $T_\mathbf{x}\mathbb{R}^2$ along the horizontal direction, and $E^u(\mathbf{x})$ be the one along the vertical direction.

**Proposition 4** *Let $\mu$ be fixed such that $0 < \mu < t_1$. For each $\mathbf{x} \in \Lambda_\mu$, there exists a positive integer $i = i(\mathbf{x})$ such that, for any $\mathbf{v} \in E^u(\mathbf{x})$*

$$\left|(d\varphi_\mu^i)_\mathbf{x}\mathbf{v}\right| > \zeta|\mathbf{v}|$$

*where $\zeta > 1$ is the constant defined above.*

**Proof.** For $0 < \mu < t_1$, $\varphi_\mu$ has a saddle fixed point $\mathbf{p}_\mu = (p_\mu, \ p_\mu)$. Let $\mathbf{x} = (x, \ y)$ be a point of $\Lambda_\mu \setminus \{\mathbf{p}_\mu\}$ satisfying $p_\mu - \delta_2 < y < p_\mu + \delta_2$. Then, there is an integer $n_0 = n_0(y)$ such that

$$n_0 = \min\{n > 0 \ : \ p_\mu + \delta_2 \leq f_\mu^n(y)\}.$$

It is clear that for each $1 \leq n < n_0$

$$\frac{d}{dy} f_\mu^n(y) > c_1$$

where $c_1 > 1$ is defined in the previous section. We here write $y_1 = f_\mu^{n_0}(y)$. From Proposition 2, there is an integer $m_1 = m_1(y_1) > 0$ such that

$$\frac{d}{dy} f_\mu^{m_1}(y_1) > \sigma_2, \quad a \leq f_\mu^{m_1}(y_1) \leq b.$$

We put $y_2 = f_\mu^{m_1}(y_1)$. From Proposition 1, there is also a finite integer $n_2 = n_2(y_2) > 0$ such that

$$\frac{d}{dy} f_\mu^{n_2}(y_2) > \sigma_1, \quad b < f_\mu^{n_2}(y_2) \leq f_\mu(\tilde{p}_\mu).$$

We rewrite $y_3 = f_\mu^{n_2}(y_2)$. From Proposition 3, there is an integer $m_3 = m_3(y_3) \geq 0$ such that

$$\frac{d}{dy} f_\mu^{m_3}(y_3) > \sigma_3, \quad l_\mu < f_\mu^{m_3}(y_3).$$

Then we have

$$\frac{d}{dy} f_\mu^{m_3+n_2+m_1+n_0}(y) > \sigma_3 \cdot \sigma_2 \cdot \sigma_1 \cdot c_1 > \sigma_3 \cdot \sigma_2 \cdot \sigma_1.$$

We write

$$i(\mathbf{x}) = i = 1 + m_3 + n_2 + m_1 + n_0.$$

For any $\mathbf{v} = (0, \ v) \in E^u(\mathbf{x})$, we have

$$\left|(d\varphi_\mu^i)_\mathbf{x} \mathbf{v}\right| = \left|\tilde{\sigma} \cdot \frac{d}{dy} f_\mu^{m_3+n_2+m_1+n_0}(y) \cdot v\right| \geq \tilde{\sigma} \cdot \sigma_1 \cdot \sigma_2 \cdot \sigma_3 \cdot |v| \geq \zeta \cdot |\mathbf{v}|.$$

For every $\mathbf{x} = (x, y) \in \Lambda_\mu$ satisfying $y \geq p_\mu + \delta_2$, by some modifications of the above estimations, we can prove it similarly. ∎

For $-t_2 < \mu < 0$, that is, after the saddle-node bifurcation, the maximal invariant set of $R_\mu$ for $\varphi_\mu$ consists of a couple of hyperbolic fixed points and a nontrivial invariant set. This nontrivial invariant set contains a saddle fixed point $\mathbf{q}_\mu = (p_\mu, q_\mu)$. For any point $\mathbf{x}$ in this nontrivial invariant set, we also take $E^u(\mathbf{x})$ and $E^s(\mathbf{x})$ which span along the vertical and horizontal directions, respectively. Since the $y$-coordinate of $\mathbf{x}$ is greater than the $y$-coordinate of $\mathbf{q}_\mu$, we can use Proposition 3 to prove the fullfillment of the expanding condition along $E^u(\mathbf{x})$. The proof is the same essentially as that of Proposition 4.

We here remember the condition of hyperbolicity: $\Lambda_\mu$ is a *hyperbolic set* for $\varphi_\mu$ if $\Lambda_\mu$ is invariant for $\varphi_\mu$ and if there is an invariant and continuous splitting $T_\mathbf{x}\mathbb{R}^2 = E^u(\mathbf{x}) \oplus E^s(\mathbf{x})$, and there exist constants $C > 0$ and $\zeta > 1$ such that, for all $\mathbf{v} \in E^u(\mathbf{x})$ and all $n > 0$, $\left|(d\varphi_\mu^n)_\mathbf{x}\mathbf{v}\right| \geq C \cdot \zeta^n \cdot |\mathbf{v}|$, and for all $\mathbf{v} \in E^s(\mathbf{x})$ and all $n > 0$, $\left|(d\varphi_\mu^n)_\mathbf{x}\mathbf{v}\right| \leq C^{-1} \cdot \zeta^{-n} \cdot |\mathbf{v}|$. See [6]. The expanding condition on $E^u(\mathbf{x})$ in the above Proposition 4 looks different from that of this definition. However, the next lemma fills this gap.

**Lemma 1** . *If for each $\mathbf{x} \in \Lambda_\mu$ there exists a finite positive integer $i = i(\mathbf{x})$ such that $\left|(d\varphi_\mu^i)_\mathbf{x}\mathbf{v}\right| > |\mathbf{v}|$ for any $\mathbf{v} \in E^u(\mathbf{x})$, then there exist constants $C > 0$ and $\zeta > 1$ such that $\left|(d\varphi_\mu^n)_\mathbf{x}\mathbf{v}\right| \geq C \cdot \zeta^n \cdot |\mathbf{v}|$ for all $\mathbf{v} \in E^u(\mathbf{x})$ and all $n > 0$.*

In [4, p.220], we find the corresponding result for the one-dimensional case. The property of two-dimensional family $\varphi_\mu$ along the vertical direction is inherited from that of the one-dimensional family $f_\mu$. Therefore we give only the outline of its proof.

**Outline of Proof of Lemma 1.** We consider the following closed sets $V_1$ and $V_2$

satisfying $\Lambda_\mu \subset V_1 \cup V_2$:

$$V_1 = \left\{(x, y) \in \mathbb{R}^2 \ : \ p_\mu - \delta_2 \leq x \leq f_\mu(\tilde{p}_\mu), \ p_\mu + \delta_2 \leq y \leq f_\mu(\tilde{p}_\mu)\right\};$$

$$V_2 = \left\{(x, y) \in \mathbb{R}^2 \ : \ p_\mu - \delta_2 \leq x \leq f_\mu(\tilde{p}_\mu), \ p_\mu - \delta_2 \leq y \leq p_\mu + \delta_2\right\}.$$

Let
$$n_1 = 1 + \min\left\{n > 0 \ : \ f_\mu^n(p_\mu + \delta_2) \geq l_\mu\right\}.$$

From the proof of Proposition 4, there is a real number $\zeta_1 > 1$ such that, for every $\mathbf{x} \in V_1$ and $\mathbf{v} \in E^u(\mathbf{x})$,

$$|(d\varphi_\mu^{n_1})_\mathbf{x} \mathbf{v}| > \zeta_1 |\mathbf{v}|.$$

Since the expanding factor of $(d\varphi_\mu)|_{V_2}$ is greater than one, there is also a real number $\zeta_2 > 1$ such that for every $\mathbf{x} \in V_2$ and $\mathbf{v} \in E^u(\mathbf{x})$,

$$|(d\varphi_\mu)_\mathbf{x} \mathbf{v}| \geq \zeta_2 |\mathbf{v}|.$$

We here let $\zeta_0 = \min\{\zeta_1, \zeta_2\}$ and $\xi = \min\{|(d\varphi_\mu)_\mathbf{x} \mathbf{v}| \ : \ \mathbf{x} \in \Lambda_\mu, \ \mathbf{v} \in E^u(\mathbf{x})\}$. Let $m$ be an integer so large that $\zeta_0^m \cdot \xi^{n_1} > 1$. We write $\tilde{n} = n_1 \cdot (m+1)$.

We here claim that $\left|(d\varphi_\mu^{\tilde{n}})_\mathbf{x} \mathbf{v}\right| > |\mathbf{v}|$ for every $\mathbf{x} \in \Lambda_\mu$ and $\mathbf{v} \in E^u(\mathbf{x})$. For any $\mathbf{x} \in \Lambda_\mu$, let $i_1 \in \{1, 2\}$ be an integer such that $\mathbf{x} \in V_{i_1}$. Moreover, let $i_2, i_3, \ldots \in \{1, 2\}$ be a sequence such that $\varphi_\mu^{n_{i_1}}(\mathbf{x}) \in V_{i_2}$, $\varphi_\mu^{n_{i_1}+n_{i_2}}(\mathbf{x}) \in V_{i_3}$ and so on. For any $\mathbf{x} \in \Lambda_\mu$, since $\tilde{n} > n_1 \cdot m$, there exist integers $s(\mathbf{x}) = s \geq m$ and $m_1 < n$ such that $\tilde{n} = n_{i_1} + n_{i_2} + \ldots + n_{i_s} + m_1$. Then, for any $\mathbf{x} \in \Lambda_\mu$ and $\mathbf{v} \in C^u(\mathbf{x})$, we have

$$\left|(d\varphi_\mu^{\tilde{n}})_\mathbf{x} \mathbf{v}\right| \geq \zeta_{i_1} \cdot \zeta_{i_2} \ldots \zeta_{i_s} \left|(d\varphi_\mu^{m_1})_{\varphi_\mu^{n_{i_1}+n_{i_2}+\ldots+n_{i_s}}(\mathbf{x})} \mathbf{v}\right| \geq \zeta_0^m \cdot \xi^{m_1} |\mathbf{v}| > |\mathbf{v}|.$$

Then, we get that $\left|(d\varphi_\mu^{\tilde{n}})_\mathbf{x} \mathbf{v}\right| > |\mathbf{v}|$ for any $\mathbf{x} \in \Lambda_\mu$. Then, there is a constant $\zeta > 1$ such that for every $\mathbf{x} \in \Lambda_\mu$, $\left|(d\varphi_\mu^{\tilde{n}})_\mathbf{x} \mathbf{v}\right| > \zeta^{\tilde{n}} |\mathbf{v}|$. For every $1 \leq i \leq \tilde{n}$, we get $\left|(d\varphi_\mu^{\tilde{n}})_\mathbf{x} \mathbf{v}\right| > C \cdot \zeta^i \cdot |\mathbf{v}|$, where $C = \min\{(\xi/\zeta)^i \ : \ 1 \leq i \leq \tilde{n}\}$. Then, for every $n > 0$, we can write $n = s\tilde{n} + t$ with $t < \tilde{n}$ and get

$$\left|(d\varphi_\mu^n)_\mathbf{x} \mathbf{v}\right| = \left|(d\varphi_\mu^{s\tilde{n}})_\mathbf{x} \cdot (d\varphi_\mu^t)_{\varphi_\mu^{s\tilde{n}}(\mathbf{x})} \mathbf{v}\right| \leq \zeta^{s\tilde{n}} \cdot C \cdot \zeta^t |\mathbf{v}| = C \cdot \zeta^n |\mathbf{v}|.$$

The property of $\varphi_\mu$ along $E^s(\mathbf{x})$ is straightforward. Any non-zero vector $\mathbf{v} \in E^s(\mathbf{x})$ is contracted linearly by $(d\varphi_\mu)_\mathbf{x}$ with a contracting factor $\lambda$ for every $-t_2 < \mu < t_1$.

From the above, we have:

**Proposition 5** . *The arc $\{\varphi_\mu\}$ has a saddle-node bifurcation for $\mu = 0$ occurring inside $\Lambda_0$. For any $-t_2 < \mu < t_1$ such that $\mu \neq 0$, $\Lambda_\mu$ is a hyperbolic set for $\varphi_\mu$,* ∎

We finally have the proof of the main theorem using the saddle-node arc $\{\varphi_\mu\}$ constructed in this section.

**Proof of Main Theorem.** Let $D_\mu = S_\mu^1 \cup R_\mu \cup S_\mu^2$. Let $\mathcal{A}$ be the $\varepsilon$-neighborhood of $\{\varphi_\mu\}$ in the $C^2$ topology. For each $\{\tilde{\varphi}_\mu\} \in \mathcal{A}$, we can write its form for any $\mathbf{x} = (x, y) \in H_\mu$ and $-t_2 < \mu < t_1$ as

$$\tilde{\varphi}_\mu(\mathbf{x}) = (\lambda(x - p_\mu) + p_\mu + \xi_1(\mu, \mathbf{x}),\ f_\mu(y) + \xi_2(\mu, \mathbf{x})),$$

where $\xi_i : \mathbb{R} \times D_\mu \to \mathbb{R}$ is a smooth map satisfying $\|\xi_i\|_2 < \varepsilon$. For every one-dimensional arc $C^2$-close to $\{f_\mu\}$, we get the same properties as are stated in Proposition 1-3, and for sufficiently small $\varepsilon > 0$, $\{\tilde{\varphi}_\mu\}$ has a saddle-node bifurcation point near 0. $\tilde{\varphi}_\mu$ is shown to be hyperbolic before and after this saddle-node bifurcation as follows. To detect a hyperbolic splitting we use the following well-known device. For $\mathbf{x} \in D_\mu$, let

$$C^u(\mathbf{x}) = \{(u, v) \in T_\mathbf{x} D_\mu : 2|u| \leq |v|\}, \quad C^s(\mathbf{x}) = \{(u, v) \in T_\mathbf{x} D_\mu : 2|v| \leq |u|\}.$$

For each $(u, v) \in C^u(\mathbf{x})$, we write $(u_1, v_1) = (d\tilde{\varphi}_\mu)_\mathbf{x}(u, v)$. For any $\varepsilon > 0$ such that $9\varepsilon < c_2 - 2\lambda$, we have

$$\frac{|u_1|}{|v_1|} = \frac{|(\lambda + \partial_x \xi_1)u + \partial_x \xi_1 v|}{|\partial_x \xi_2 u + \partial_y(f_\mu + \xi_2)v|} < \frac{(\lambda + \epsilon)|u| + \varepsilon|v|}{(c_2 - \varepsilon)|v| - \varepsilon|u|} < \frac{1}{2},$$

that is, $(d\tilde{\varphi}_\mu)_{\mathbf{x}} C^u(\mathbf{x}) \subset C^u(\tilde{\varphi}_\mu(\mathbf{x}))$. For each $\mathbf{v} = (u, v) \in C^u(\mathbf{x})$, we use the norm $\|\mathbf{v}\| = \max\{|u|, |v|\}$. Since $\tilde{\varphi}_\mu$ is $C^2$-close to $\varphi_\mu$, we get $\|(d\tilde{\varphi}_\mu^i)\mathbf{v}\| > \zeta\|\mathbf{v}\|$. Similarly, we have $(d\tilde{\varphi}_\mu)_{\mathbf{x}}^{-1} C^s(\mathbf{x}) \subset C^s(\tilde{\varphi}_\mu^{-1}(\mathbf{x}))$, and $\|(d\tilde{\varphi}_\mu)\mathbf{v}\| < \zeta^{-1}\|\mathbf{v}\|$ for any $\mathbf{v} \in C^s(\mathbf{x})$. For any $\mathbf{x} \in \widetilde{H}$, the proof is even simpler. ∎

## Acknowledgments

We would like to thank Professor Atsuro Sannami and Professor Floris Takens for their many valuable suggestions for improvements of this paper. This work was partially supported by NNSF of China, and the Reserch Institute for Technology of Tokyo Denki University under contract no. Q97J-05.